\documentclass[11pt]{amsart}
\usepackage{amssymb,amsmath,epsf}

\input xy
\xyoption{all}

\newtheorem{thm}{Theorem}[section]
\newtheorem{lem}[thm]{Lemma}

\newtheorem{prop}[thm]{Proposition}
\newtheorem{exe}[thm]{Example}

\theoremstyle{definition}
\newtheorem{defn}[thm]{Definition}

\def\bZ{\mathbb{Z}}

\def\bR{\mathbb{R}}

\newcommand{\ra}{\rightarrow}

\newcommand{\bdy}{\partial}

\newcommand{\Mforms}[1]{\Omega^{#1}M}
\newcommand{\iMforms}[1]{\Omega_\Lambda^{#1}M}

\newcommand{\rizaMforms}[1]{\widetilde{\Omega}^{#1} M}

\newcommand{\Aforms}[1]{\Omega^{#1}A}

\newcommand{\iAforms}[1]{\Omega^{#1}_\Lambda A}
\newcommand{\imAforms}[1]{\Omega_{Im}^{#1}A}

\newcommand{\MAforms}[1]{\Omega^{#1}(\rho)}

\newcommand{\iMAforms}[1]{\Omega_\Lambda^{#1}(\rho)}
\newcommand{\reduct}[1]{\widetilde {#1}}

\newcommand{\MAint}[1]{C^{#1}(\rho;\Lambda)}
\newcommand{\MAintcoh}[1]{H^{#1}(C_\rho;\Lambda)}

\newcommand{\Mreal}[1]{C^{#1}(M;\bR)}
\newcommand{\Mrealmodint}[1]{C^{#1}(M;\bR / \Lambda)}
\newcommand{\MAreal}[1]{C^{#1}(\rho;\bR)}
\newcommand{\MArealcoh}[1]{H^{#1}(C_\rho;\bR)}

\newcommand{\csg}[2]{\widehat{H}^{#1} ({#2})}
\newcommand{\MAcsg}[1]{\widehat{H}^{#1}(\rho)}

\newcommand{\MAintcsg}[1]{\widehat{H}_0^{#1}(\rho)}

\newcommand{\MAchat}[1]{\widehat{C}^{#1}(\rho)}
\newcommand{\iMAchat}[1]{\widehat{C}_0^{#1}(\rho)}
\newcommand{\MAzhat}[1]{\widehat{Z}^{#1}(\rho)}

\newcommand{\MAccheck}[1]{\check{C}^{#1}(\rho)}

\newcommand{\Mzcheck}[1]{\check{Z}^{#1}(M)}

\newcommand{\MAhsg}[1]{\check{H}^{#1}(\rho)}

\newcommand{\ccheck}[2]{\check{C}^{#1} ({#2})}

\def\cocycle{((c,b),(h,e),(\omega,\theta))}

\DeclareMathOperator{\Hom}{Hom}

\DeclareMathOperator{\Ker}{Ker}

\DeclareMathOperator{\Image}{Im}

\def\ra{\rightarrow}

\title{Relative differential characters}
\author{Mark Brightwell} \address{129 rue de Lausanne, 1202 Geneva,
  Switzerland} 
\author{Paul Turner}
\address{School of Mathematical and Computer
  Sciences \\Heriot-Watt University\\ Edinburgh EH14 4AS\\Scotland}

\begin{document}


\begin{abstract}
There are two natural candidates for the group of relative
Cheeger-Simons differential characters. The first directly extends the
work of Cheeger and Simons and the second extends the description
given by Hopkins and Singer of the
Cheeger-Simons group as the homology of a certain cochain complex. 
We discuss both approaches and relate the two relative groups.
\end{abstract}

\maketitle


\section{Introduction}
Given a principal $G$-bundle with connection on a smooth manifold $M$,
an invariant polynomial gives rise to a Chern-Weil form in the
base. Lifting this to the total space gives an exact form and
the well known form of Chern and Simons \cite{ChernSimons} has the
property that its differential is this exact form. If one happens to have a global
section one can pull this back to the base space, but in general
Chern-Simons forms live in the total space. Cheeger and Simons
\cite{CheegerSimons} introduced the notion of a differential character as a way of
understanding Chern-Simons forms in terms of the base space. These differential
characters can be assembled into a
group $\csg k M$ which contains both homotopical and geometrical
information. One may also think of differential characters as encoding
how integral cycles and differential forms with integral periods
interact when viewed in real cohomology. When $k=2$ the Cheeger-Simons
group classifies $U(1)$-bundles with connection and when $k=3$ it
classifies $U(1)$-gerbes with connection (see for example \cite{Brylinski}).
The latter makes Cheeger-Simons groups relevant to ``stringy''
topology since a gerbe with connection can be interpreted as a line
bundles over a free loop space with parallel transport over surfaces. In
order to incorporate open strings into such a framework it is
necessary to develop  a relative theory.

In this paper we study the notion of {\em relative} differential
character associated to a pair of manifolds $A$ and $M$ together with
a smooth map $\rho\colon A \ra M$. There are two competing definitions
for the relative group: the first a natural extension of Cheeger and
Simons' work giving what we will call the relative Cheeger-Simons
group, the second following the more recent work of Hopkins and Singer
defining the Cheeger-Simons group as the cohomology of a certain
cochain complex, giving what we will call the relative Hopkins-Singer
group. The purpose of this paper is to clarify the relationship
between these two groups.

For compact smooth manifolds $A$ and $M$ and a smooth map $\rho:A\ra
M$, we give in Section \ref{section:rel-csg} a natural definition of
the group of relative differential characters $\MAcsg k$, which
closely mirrors the original definition of Cheeger and Simons.  We
show that this group can be described as the homology of a cochain
complex analogous to the description of (non-relative) Cheeger-Simons
characters by Hopkins and Singer \cite{HopkinsSinger}. The relative
group fits into certain short exact sequences as in the non-relative
case. We end section 2 with some examples. There is, on the other
hand, another natural candidate for a relative group if one adopts the
Hopkins-Singer point of view from the outset. In Section 3 we define
the relative Hopkins-Singer group $\MAhsg {k}$.  Unlike the relative
Cheeger-Simons group this group fits into a long exact sequence.  Finally,
in Section 4 we relate the relative Cheeger-Simons group to the
relative Hopkins-Singer group showing that the latter is a quotient of
a certain subgroup of the former.


\section{Relative Cheeger-Simons groups}\label{section:rel-csg}
We begin by briefly reviewing the definition of differential
characters given in \cite{CheegerSimons} and provide a variant of the
Hopkins-Singer description as the homology of a certain cochain
complex \cite{HopkinsSinger}. We then introduce the group of relative
differential characters, provide a description as the homology
of a complex and present three short exact sequences. Finally in
this section we give some examples of relative differential
characters. We now fix for the rest of the paper a proper subring
$\Lambda$ of $\bR$.

\subsection*{Cheeger-Simons differential characters}
 Let $M$ be a compact closed smooth manifold. We denote by $C_kM$ and
$Z_kM$ the groups of smooth $k$-chains and $k$-cycles in $M$. We let
$\Mforms {k}$ denote the group of $k$-forms on $M$ which can be viewed
as a subgroup of $\Mreal {k}$ via integration. 
We let $\iMforms
{k}$ denote the subgroup of $k$-forms with
$\Lambda$-periods, noting also that such forms are closed.
One important property that we shall use many times is 
that a non-zero differential form can
never take values only in a proper subring. 
It follows that  we can view $\Mforms {k}$
as a subgroup of $\Mrealmodint {k}$ and we will make this
identification implicitly throughout.

The group of Cheeger-Simons differential characters
of degree $k\geq 1$ is defined as
\[
\csg k M = \{ f\in \Hom (Z_{k-1}M,  \bR / \Lambda)  \mbox{ such that } f\circ \partial \in
\Mforms k\}.
\]
Cheeger and Simons index this group by $k-1$, but it is convenient to
shift up the dimension by 1 to be consistent with Hopkins and
Singer. The example one should have in mind is when $k=2$, where a
line bundle with connection gives a differential character via its holonomy.

These groups sit in three exact sequences as follows. 
\begin{eqnarray*}
&&0 \ra H^{k-1}(M; \bR / \Lambda) \ra \csg k M   \ra  \iMforms k \ra 0\\
&&0 \ra {\Mforms {k-1}}/{\iMforms {k-1} } \ra \csg k M \ra H^k(M;\Lambda) \ra 0\\
&&0 \ra {H^{k-1}(M;\bR)}/{r(H^{k-1}(M;\Lambda))} \ra \csg k M \ra R^k(M;\Lambda) \ra 0,
\end{eqnarray*}
where $r$ is the natural map $H^{k-1}(M;\Lambda) \ra H^{k-1}(M;\bR)$ and 
\[
R^k(M;\Lambda)=\{(\omega ,u)\in \iMforms k \times H^k(M;\Lambda) \mid [\omega]=r(u)\}.
\]

The closed differential form associated to the differential character
is to be thought of as a curvature form and the class in
$H^k(M;\Lambda) $ as a characteristic class analogous to the first
Chern-Class for $U(1)$-bundles.

We now turn to the description of Cheeger-Simons groups as the
homology of a cochain complex following Hopkins and Singer.
Consider the cochain complex $\ccheck * M $ defined by 
\[
\ccheck k M = C^k(M,\Lambda) \times C^{k-1}(M;\bR) \times \iMforms k
\]
with differential $\Delta \colon \ccheck k M \ra \ccheck {k+1} M$ given by
\[
\Delta(c,h,\omega) = (\delta c, \omega - c - \delta h , 0).
\]
This differs from the Hopkins-Singer approach in that we include only
forms with $\Lambda$-periods in the complex above. This allows us to
obtain all Cheeger-Simons groups as the homology of a single chain
complex. Indeed, denoting the subgroup of cocycles by $\Mzcheck k$, we can
define a map $\varphi \colon \Mzcheck k \ra \csg {k}{M}$  by   
 $\varphi (c,h,\omega)=\widetilde {h|}_{Z_{k-1}M}$ (here  the tilde means
mod $\Lambda$ reduction) and  Hopkins and Singer show that $\phi$  induces an isomorphism
\begin{equation*}
\label{thm:HS-CS-iso}
H_k(\ccheck * M) \simeq \csg {k}{M}.
\end{equation*}

\subsection*{Relative Cheeger-Simons groups}
Let $A$ and $M$ be compact smooth manifolds, possibly with boundary, and
$\rho:A\ra M$ a smooth map. In order to define a relative version of Cheeger-Simons
we find it convenient to consider the following chain complex:
\[
C_k(\rho)=C_kM \times C_{k-1}A 
\]
with differential $\bdy : C_k(\rho)\ra C_{k-1}(\rho)$ given by
\[
\bdy (\sigma, \tau)=(\bdy \sigma +\rho_*\tau ,-\bdy \tau).
\]
We will denote the cycles and boundaries of this complex by
$Z_k(\rho)$ and $B_k(\rho)$ respectively. For a general map $\rho$ the
homology of this complex is the homology $H_k(C_\rho)$ of the mapping
cone $C_\rho$, while in the special case of submanifold $A\subset M$ with
$\rho$ the inclusion map we recover the usual relative homology groups
$H_k(M,A)$.

The associated dual complex $C^*(\rho;G) = \Hom (C_*(\rho),G)$ with
coefficient group $G$ has differential $\delta (h,e)=(\delta
h,\rho^*h-\delta e)$.  In general, taking homology gives
$H^k(C_\rho;G)$ and when $\rho$ is an inclusion this is isomorphic to
$H^k(M,A;G)$.  Notice the symbols $\bdy$ and $\delta$ are used both to
denote the relative differentials and as the usual operators in the
non-relative complexes.

We will also use the complex of pairs of differential forms
\[
\MAforms * =\Mforms * \times \Aforms {*-1}
\]
with differential $\delta(\omega, \theta)=(d\omega, \rho^*\omega
-d\theta)$. The homology of this
complex gives the relative de Rham cohomology (see for instance
\cite{BottTu}). As above in the non-relative version we can view
 $\MAforms k$ as a subgroup of $ C^k(\rho;\bR/\Lambda)$ via integration: for $(\omega,\theta)\in
 \MAforms k$ set
\[
(\omega,\theta)(\sigma,\tau)=\int_\sigma \omega + \int_\tau \theta
\]
where $(\sigma,\tau)\in C_k(\rho)$. For $\rho$ the inclusion map, this induces the relative de
Rham isomorphism. 

There is now an obvious definition to be made for the notion of {\em
relative} differential character.

\begin{defn} 
\label{defn:rel-csg}
Let $\rho:A\ra M$ be a smooth map. The group of {\em relative Cheeger-Simons differential characters} of degree
$k$ associated to the triple $(A,\rho , M)$ is defined as
\[
\MAcsg k=\{f\in \Hom (Z_{k-1}(\rho),  \bR/\Lambda) \text{ such that } f\circ \bdy \in \MAforms k \}.
\]
\end{defn}
 This definition has already appeared in \cite{Zucchini} in the case
 $\rho$ is an inclusion and a close
variant has also appeared in \cite{HarveyLawson} for the
restricted case of manifolds with boundary. It is immediately seen
that if we take $A=\emptyset$ we recover the non-relative group $\csg
{k}{M}$ and we present two further examples at the end of this
section.

Associated to a differential character is a pair of differential forms
and a characteristic class in the cohomology of
the mapping cone. We shall denote by $\iMAforms {k+1} <\MAforms {k+1}$ the subgroup of
pairs $(\omega,\theta)$ taking $\Lambda$ values on relative cycles and say that
such pairs have {\em relative $\Lambda$-periods}. In fact if $\rho$ is the inclusion, $\iMAforms
{*}$ consists of precisely  those forms representing the image of
$\Lambda$-classes in relative real cohomology (via the de Rham
isomorphism). In particular, a pair with relative $\Lambda$-periods satisfies the following properties.

\begin{lem}\label{lem:relperiods}
If $(\omega,\theta)$ has relative $\Lambda$-periods then 
\begin{enumerate}
\item $\delta (\omega,\theta)=0$ in the complex  $\MAforms *  $
  i.e. $d\omega = 0$ and $d\theta = \rho^*\omega$.
\item $\omega$ has $\Lambda$-periods.
\end{enumerate}
\end{lem}

\begin{proof}
For (1) let $\alpha \in C_{k}A$. Then
\[
\int_\alpha (\rho^*\omega-d\theta) = \int_{\rho_*\alpha}\omega  - \int_{\partial \alpha} \theta= (\omega, \theta)(\rho_*\alpha, -\partial \alpha)
\] 
Since $(\omega, \theta)$ has relative
$\Lambda$-periods and $(\rho_*\alpha,-\partial \alpha)$ is a relative cycle, the right hand side is in $\Lambda$. Recalling 
that a non-zero form never takes values only in a proper subring we
conclude that $d\theta - \rho^*\omega= 0 $. Showing $d\omega = 0$ is
similar.

For (2) let $\sigma\in Z_kM$ then
\[
\int_\sigma \omega =\int_\sigma \omega - \int_0\theta =
(\omega,\theta)(\sigma,0) 
\]
and the right-hand side in an element of $\Lambda$ since
$(\omega,\theta)$ has relative $\Lambda$-periods.
\end{proof}

We now define two homomorphisms
\begin{eqnarray*}
&&\delta_1 :\MAcsg {k} \rightarrow \iMAforms {k},\\
&&\delta_2 :\MAcsg {k} \rightarrow \MAintcoh {k}.
\end{eqnarray*}
Let $f\in \MAcsg k$. Since $\bR/\Lambda$ is divisible and
$Z_{k-1}(\rho)$ free, $ f$ can be extended to a map
$C_{k-1}(\rho)=C_{k-1}M \times C_{k-2}A \rightarrow \bR$, which we
denote by the pair $(h,e)$. From  the definition of differential
character there exists a pair $(\omega ,\theta)\in \MAforms k$ and and
pair $(c,b)\in
\MAint k$ such that
\[
\delta(h,e)=(\omega,\theta)-(c,b).
\]
The above homomorphisms are defined by
\[ 
\delta_1(f)=(\omega,\theta) \;\;\;\;\;\;\;\;\;\; \delta_2(f)=[c,b].
\]
To see that $\delta_1$ is well defined observe that for a  relative
cycle $(\sigma,\tau)\in Z_k(\rho)$ we have,
\[
(\omega,\theta)(\sigma,\tau)=(\delta(h,e) + (c,b))(\sigma,
\tau)= (c,b)(\sigma,\tau)\in \Lambda
\]
and hence $(\omega,\theta)$ has relative $\Lambda$-periods. For $\delta_2$ observe that 
\[
\delta(c,b) = \delta((\omega, \theta) - \delta(h,e)) = \delta(\omega,
\theta) - 0 = 0
\]
since $(\omega,\theta)$ is closed by Lemma \ref{lem:relperiods}. Moreover
$\delta_1$ and $\delta_2$ are independent of the choice of the
extension $(h,e)$ of $f$.

It would  now
possible to determine the kernels of $\delta_1$ and $\delta_2$ 
directly and fit
them into short exact sequences. We prefer to first describe the relative
Cheeger-Simons groups as the homology of a cochain complex along the
lines of the non-relative
case described above. We can then easily deduce the
short exact sequences using this description.

Consider the complex
\[
\MAchat {*} =\MAint {*} \times \MAreal {*-1} \times \iMAforms {*}
\]
with differential
\[
\Delta((c,b),(h,e),(\omega,\theta))=(\delta(c,b),(\omega,\theta)-(c,b)-\delta(h,e),0).
\]

The last entry is in fact $\delta (\omega ,\theta)$ but this is zero since
$(\omega,\theta)\in \iMAforms *$ and hence is closed. Let $\MAzhat k$ denote
the group of $k$-cocycles in $\MAchat {*} $ and define a map 
\[
\varphi \colon \MAzhat k \ra \MAcsg
{k}
\]
 by
\[
\varphi ((c,b),(h,e),(\omega,\theta))=\reduct {(h,e)}|_{Z_{k-1}(\rho)}
\]
where as above the tilde denotes mod $\Lambda$ reduction.  Since we
start with a cocycle we have $\delta(h,e)=(\omega,\theta)-(c,b)$ and
so $\widetilde{(h,e)}\circ \bdy=\widetilde{(\omega,\theta)}$ which
shows $\widetilde{(h,e)}|_{Z_{k-1}(\rho)}$ is indeed a differential
character. 

\begin{thm}
\label{thm:rel-CS-HS-iso}
The homomorphism $\varphi$ defined above induces an isomorphism
\[ H_k(\MAchat {*}) \cong \MAcsg {k}. \]
\end{thm}
\begin{proof}
To show $\varphi$ is surjective let $f:Z_{k-1}(\rho)\ra \bR/\Lambda$
be a relative differential character. From the discussion of the maps
$\delta_1$ and $\delta_2$ above, we can find pairs $(h,e)\in \MAreal
{k-1}$, $(c,b)\in C^k(\rho;\Lambda)$ and $(\omega,\theta)\in \iMAforms
k$ such that $\delta(c,b) = 0 $ and
$\delta(h,e)=(\omega,\theta)-(c,b)$. In other words $\cocycle$ is a
cocycle which maps to $f$ under $\varphi$.

The proof will be completed by showing that $ \Ker \varphi = \Image
\Delta$. 

Firstly, an element in the image of $\Delta$ has the form
$(\delta(c^\prime,b^\prime),(\omega^\prime,\theta^\prime)-(c^\prime,b^\prime)-
\delta(h^\prime,e^\prime),0)$ where
$((c^\prime,b^\prime),(h^\prime,e^\prime),
(\omega^\prime,\theta^\prime))\in \MAchat {k-1}$. Since
$(\omega^\prime,\theta^\prime)$ has relative $\Lambda$-periods and
$(c^\prime, b^\prime)$ is $\Lambda$-valued, this is clearly
mapped to 0 under $\varphi$. Hence $\Image \Delta \subset \Ker \varphi$.

Now suppose $\cocycle \in \Ker \varphi$. Then
$\reduct{(h,e)}|_{Z_{k-1}(M,A)}=0$, so $\delta\reduct
{(h,e)} = 0$. Since $(\omega,\theta) = (c,b) + \delta (h,e)$ we have
\[
\reduct{(\omega,\theta)}= \reduct{(c,b)}+ \delta\reduct{(h,e)} = 0
\]
as so $(\omega,\theta)=0$. Using the splitting of the standard short
exact sequence relating cycles, chains and boundaries we can extend
the $\Lambda$-valued map $-(h,e)|_{Z_{k-1}(\rho)}\colon Z_{k-1}(\rho)
\ra \Lambda$ to a map $C_{k-1}(\rho)\ra \Lambda$. Denoting this
extension by $(c^\prime, b^\prime)$ we have
\[
\delta(c^\prime,b^\prime) = -\delta(h,e)= (c,b) -
(\omega, \theta)= (c,b).
\] 
Finally, since
$-(c',b')-(h,e)$ is zero on cycles it is a real coboundary, namely
$-(c',b')-(h,e)=\delta (h',e')$ for some $(h',e')\in \MAreal
{k-2}$. Hence 
\[
\Delta((c',b'),(h',e'),0)=
(\delta(c^\prime,b^\prime),-(c^\prime,b^\prime)-\delta(h',e'),0 )=
((c,b),(h,e),(\omega,\theta)),
\]
and so $\cocycle \in \Image \Delta$. Thus $\Ker \varphi \subset
\Image \Delta$ which finishes the proof.
\end{proof}

Using this description is is now easy to fit the relative
Cheeger-Simons groups into three short exact sequences. Let
$r:\MAintcoh {k}\ra \MArealcoh {k}$ be the natural map induced by the
inclusion $\Lambda \ra \bR$ and set 
\[
R^k(\rho,\Lambda)=\{ ((\omega,\theta),u)\in \iMAforms k \times
\MAintcoh k \mid [\omega,\theta]=r(u) \}
\]
where  $[\omega,\theta] $ denotes the de Rham cohomology class.

\begin{thm}
\label{thm:csg_ses}
The following sequences are exact.
\[
\xymatrix{0\ar[r] &  H^{k-1}(C_\rho;\bR/\Lambda)\ar[r]& \MAcsg {k}
  \ar[r]^{\delta_1}& \iMAforms {k} \ar[r] &  0\\
 0\ar[r]&  \MAforms {k-1}/\iMAforms {k-1} \ar[r]& \MAcsg {k}
  \ar[r]^{\delta_2}& \MAintcoh {k} \ar[r] & 0\\
 0\ar[r]& \MArealcoh {k-1}/r(\MAintcoh {k-1}) \ar[r]& \MAcsg {k}
  \ar[r]^{(\delta_1,\delta_2)}&  R^{k}(\rho;\Lambda)\ar[r]& 0
}
\]
\end{thm}

\begin{proof}
Consider the following short exact sequences of chain complexes.
\begin{eqnarray*}
&&0\ra \MAint {*} \times \MAreal {*-1} \ra \MAchat {*} \ra \iMAforms {*} \ra 0\\
&&0\ra \MAreal {*-1} \times \iMAforms {*} \ra \MAchat {*} \ra \MAint {*} \ra 0\\
&&0\ra \MAreal {*-1} \ra \MAchat {*} \ra \MAint {*} \times \iMAforms {*} \ra 0
\end{eqnarray*}
Each of the differentials in the complexes on the  left is  obtained by
restriction of the differential of $\MAchat {*}$. Then as in
\cite{HopkinsSinger} the three exact sequences of the theorem can be
deduced from the long exact sequences obtained from the short exact
sequences above.
\end{proof}

Notice that taking $A = \emptyset$ recovers the usual non-relative
short exact sequences. Unfortunately, there is no long exact sequence
for the relative group of the triple $(A,\rho, M)$. This is easily
seen just at the level of differential forms: a differential character
with associated pair of forms $(\omega,\theta)$ induces a differential
character on $M$ with form $\omega$, however unless $\theta$  is
closed, there is no guarantee that $\rho^*
\omega = 0$ which would be required by any long exact sequence.

\subsection*{Examples}
We end this section by presenting a couple of examples. In both 
$\rho$ is the inclusion map and we adopt the standard relative
notation $Z_k(M,A)$ and $\widehat{H}^k (M,A)$ for
$Z_k(\rho)$ and $\MAcsg k$.

\begin{exe}
Let $M$ be a smooth manifold with submanifold $A$. Let $P\ra M$ be a
principal $S^1$-bundle with connection together with a given trivialisation
of the bundle over $A$. The relative holonomy
$H:Z_1(M,A)\ra S^1$ of $P$ is defined as follows. Given $\gamma\in Z_1(M,A)$, split $\gamma$
into a (finite) collection of piecewise smooth paths $\gamma_i$ with
endpoints in $A$. For each $\gamma_i$, define $H(\gamma_i)$ to be the
holonomy along $\gamma_i$ where the endpoint fibres are identified via
the given trivialisation  over $A$. Then $H(\gamma)$ is defined to be the
product (in the circle group) of the $H(\gamma_i)$.

We now define the relative differential character $\widehat\chi:
Z_1(M,A)\ra \bR/\bZ$ associated to $P$ by the equation
\[
e^{2\pi i \widehat\chi(\gamma)}=H(\gamma).
\] 
To see that $\widehat\chi$ is indeed a differential character consider
the curvature 2-form $\omega$ of $P$. Since $P$ is trivial over $A$
and hence $\omega$ exact on $A$ there exists a 1-form $\theta$ such
that on $A$ we have $\omega=d\theta$. This form satisfies
$H(\tau)=e^{2 \pi i  \int_\tau \theta}$ for $\tau \in C_1A$, which in
turn immediately implies that $H(\bdy
(\sigma,\tau))=e^{2\pi i  (\int_\sigma \omega +\int_\tau\theta)}$ for a relative boundary
$\bdy(\sigma,\tau)$. Using the defining equation for $\widehat\chi$,
we see that
\[
\widehat\chi (\bdy(\sigma,\tau))=\int_\sigma \omega
  +\int_\tau\theta =(\omega,\theta)(\sigma,\tau) \mod \bZ
\]     
and hence $\widehat\chi \in \widehat H^2(M,A)$.
\end{exe}

\begin{exe}
In this example we consider $A=S^{n-1}$ as the boundary of
$M=D^n$. Using the exact sequences in Theorem \ref{thm:csg_ses}
and some elementary topology one obtains the following.
\[
\widehat H^k(D^n,S^{n-1}) = 
\begin{cases}
C^\infty (D^n, \bR/\Lambda) & k=1 \\
\Omega^k_\Lambda (D^n,S^{n-1}) & 2\leq k \leq n\\
\bR / \Lambda & k=n+1 \\
0 & k\geq n+2
\end{cases}
\]

\end{exe}


\section{Relative Hopkins-Singer groups}
Following the description of Hopkins and
Singer there is another candidate for the relative group. This group, which we refer to
as the relative Hopkins-Singer group can be identified
with  a quotient of a subgroup of the
relative Cheeger-Simons group and does fit into a 
long exact sequence. 

Recall from the previous section the chain complex
\[
\ccheck k M = C^k(M,\Lambda) \times C^{k-1}(M;\bR) \times \iMforms k
\]
with differential $\delta (c,h,\omega)=(\delta c,\omega -c-\delta
h,0)$ which, as we have seen,  has homology $\csg {*}{M}$. Following
the approach used to define relative real and de Rham cohomology
consider the chain complex
\[
\ccheck * M \times \ccheck {*-1} A
\]
with differential $\delta (S,T)=(\Delta S, \rho^*S-\Delta T)$, where
$S,T$ are triples in $\ccheck * M$ and $\ccheck {*-1} A$
respectively. 
\begin{defn} The \emph{relative Hopkins-Singer
groups} are defined by
\[\MAhsg {k}= H_k(\ccheck * M \times \ccheck {*-1} A)\]
\end{defn}

These groups can be fitted into short exact sequences analogous to
those above for the relative Cheeger-Simons groups. However, the
groups of forms appearing are not particularly enlightening and it is
better to understand this group in terms of its relation to the relative
Cheeger-Simons group which is the subject of the next section.
  
A main feature of the relative Hopkins-Singer groups is that they do
fit into a long exact sequence. Let $q$ be the projection $\ccheck * M
\times \ccheck {*-1} A \ra \ccheck * M$ and let $l$ be the inclusion
$\ccheck {*-1} A \ra \ccheck * M \times \ccheck {*-1} A$. Using the
same notation for induced maps we have the following theorem.

\begin{thm}
The following sequence is exact.
\[
\xymatrix{ \dots \ar[r] & \csg {k-1}{A} \ar[r]^{l} & \MAhsg
{k}\ar[r]^{q} & \csg {k}{M} \ar[r]^{\rho} & \csg {k}{A}\ar[r] & \dots}
\]
where $\rho$ is the natural pullback map from differential
characters on $M$ to those on $A$.
\end{thm}
\begin{proof}
This is the long exact sequence induced by the following exact sequence of chain complexes
\[
0\ra \ccheck {*-1} A \ra \ccheck * M \times \ccheck {*-1} A \ra \ccheck * M \ra 0
\]
with the boundary map of this sequence coinciding with the map $\rho$ above.
\end{proof}


\section{Relating the relative Hopkins-Singer groups and the relative Cheeger-Simons groups}

We now relate the relative Hopkins-Singer group to the relative
Cheeger-Simons group by
identifying $\MAhsg k$ as a quotient of a subgroup of $\MAcsg k$.

In order to do this it
will be helpful to view the relative Hopkins-Singer groups as the
homology of a slightly different chain complex. Consider the chain complex
\[
\MAccheck {*}=\MAint {*}\times \MAreal {*-1}\times \iMforms {*}\times \iAforms {*-1}
\]
with differential $\Delta :\MAccheck {*}\ra \MAccheck {*+1}$ given by
\[
\Delta((c,b),(h,e),(\omega,\theta))=(\delta(c,b),(\omega,\theta)-(c,b)-\delta(h,e),(0,\rho^*\omega))
\]
Notice that the last term is in fact $\delta (\omega,\theta)=(d \omega
,\rho^*\omega -d\theta)=(0,\rho^*\omega)$.

This chain complex $\MAccheck {*}$ is isomorphic (by
changing the sign of $e$) to $\ccheck *
M \times \ccheck {*-1} A$ and so we have
\[
\MAhsg {k}\simeq H_k(\MAccheck {*}).
\]

 Let
$\MAintcsg k< \MAcsg k$ be the subgroup of relative differential
characters defined by
\[
\MAintcsg k = \{f\in \MAcsg k \mbox{ such that } f|_{Z_{k-1}A}=0\}.
\]
Here we are regarding $Z_{k-1}A \subset Z_{k-1}M \subset
Z_{k-1}(M,A)$ i.e. regard $\tau\in Z_{k-1}A$ as $(\tau,0)\in
Z_{k-1}(M,A)$. Notice that 
\[
f(\tau,0) = f\circ \partial (0,\tau) =
(\omega, \theta)(0,\tau) = \widetilde{\int_\tau \theta},
\]
so the condition
$f|_{Z_{k-1}A}=0 $ is equivalent to $\theta$ having $\Lambda$-periods.

These groups can also be identified with the cohomology of a cochain
complex. Let
\[
\rizaMforms k = \{ \omega \in \Mforms k \mid (\omega,0)\in \iMAforms k \}.
\]
noticing that $ \rizaMforms k \subset \iMforms {k}$. Then $\MAintcsg
k$ is the homology of the cochain complex
\[
\iMAchat * = \MAint {*}\times \MAreal {*-1}\times
\rizaMforms * \times \iAforms {*-1}
\]
with differential 
\[
\Delta((c,b),(h,e),(\omega,\theta))=(\delta(c,b),(\omega,\theta)-(c,b)-\delta(h,e),0).
\]

We can easily relate $\MAintcsg k$ to  $\MAcsg k$. Let
\[
\imAforms {k-1} = \Image ( \iMAforms k \ra \Aforms {k-1}).
\]
and note that given $\theta\in \iAforms {k-1}$ then $(0,\theta)\in
\iMAforms k $ and so $\iAforms {k-1} < \imAforms {k-1}$.

\begin{prop}\label{prop:intrelcs}
The following sequence is exact.
\begin{equation*}
0\ra \MAintcsg k\ra  \MAcsg k  \ra \imAforms {k-1}/\iAforms {k-1} \ra 0
\end{equation*}
where the left map is the obvious inclusion and the right map the
projection.
\end{prop}

\begin{proof}
Consider the short exact sequence of complexes
\[
\xymatrix{ 0 \ar[r] & \iMAchat * \ar[r]^{j} & \MAchat * \ar[r]^-{p} &
  \imAforms {*-1}/\iAforms {*-1} \ar[r] & 0}
\]
where complex on the right has trivial differential and $j$ and $p$
are the inclusion and projection maps. To show this is exact the
only point worth noting is that given $\theta\in\imAforms {k-1}$ then
by definition there exists $\omega$ such that $(\omega,\theta)\in
\iMAforms k$ and so $((0,0),(0,0),(\omega,\theta))$ is a cocycle mapped to
$\theta$ under $p$ i.e. $p$ is onto. 

The differential of the associated long exact sequence is trivial
giving the result.
\end{proof}

To relate $ \MAintcsg k$ to the relative Hopkins-Singer
group we need to define two maps $\phi$ and $J$. We define the  map
\[
\phi :\iMforms {k-1}/ \rizaMforms {k-1} \ra \MAintcsg k
\]
as follows. For $\nu \in \iMforms {k-1}$ define $\phi(\nu) :Z_{k-1}(\rho)\ra
\bR/\Lambda$ by
\[
\phi(\nu)(\gamma,-\bdy \gamma)=\int_\gamma
\nu \;\;\;\;\;\mbox{ mod $\Lambda$ }.
\]
If $(\sigma,\tau)\in C_k M\times C_{k-1} A$ we have $\phi(\nu) \circ
\bdy (\sigma,\tau)= (0,\rho^*\nu)(\sigma,\tau)$ and so $\nu(\rho)\circ
\partial \in \MAforms k$. Thus $\phi(\nu)$ is a differential
character as required and moreover if $\nu\in \rizaMforms {k-1}$ then
$\rho^*\nu = 0$ and so $\phi$ is well defined.

To define $J$, recall from Section 2 that given a relative
differential character $f$ we can find pairs $(c,b),(h,e)$ and
$(\omega,\theta)$ such that $((c,b),(h,e),(\omega,\theta))$ is a
cocycle in $\MAchat *$. If $f\in \MAintcsg k$ then in addition we have
$\theta$ has $\Lambda$-periods and $\rho^*\omega = 0$. Since $\omega$
also has $\Lambda$-periods (by Lemma \ref{lem:relperiods}) it follows
that $((c,b),(h,e),(\omega,\theta))$ is in fact a cocycle in
$\MAccheck *$. Thus we can define a map
\[
J\colon \MAintcsg k  \ra \MAhsg k 
\]
by 
\[
J(f) = [((c,b),(h,e),(\omega,\theta))].
\]

The relationship between the relative Cheeger-Simons groups and the
relative Hopkins-Singer groups is given by Proposition
\ref{prop:intrelcs} and the following result.

\begin{thm}
The following sequence is exact.
\[
\xymatrix{0 \ar[r] &  \iMforms {k-1}/ \rizaMforms {k-1}
  \ar[r]^-{\phi} &  \MAintcsg k  \ar[r]^{J} & \MAhsg k \ar[r] & 0}
\]
\end{thm}

\begin{proof}
Consider the short exact sequence of complexes
\[
\xymatrix{ 0 \ar[r] & \iMAchat * \ar[r]^{j} & \MAccheck * \ar[r]^-{p} &
  \iMforms {k}/ \rizaMforms {k} \ar[r] & 0}
\]
where the complex on the right has trivial differential, $p$ is the
projection and $j$ is the inclusion. Note that $j$ is well defined since by
Lemma \ref{lem:relperiods} $\omega \in \iMforms k $ and by definition
$\theta\in \iAforms {k-1}$. It is also a map of cochain complexes
since for $((c,b),(h,e),(\omega,\theta))\in \iMAchat k$ we have
$\rho^*\omega = d\theta = 0$.

The only non-trivial part in showing that the above sequence is indeed
exact is the middle term. If
$((c,b),(h,e),(\omega,\theta))\in \Ker (p)$ then $\omega\in
\rizaMforms {k}$ i.e. $(\omega, 0 ) \in \iMAforms k$. Note that since
$\theta$ has $\Lambda$-periods we have $(0,\theta)\in\iMAforms
k$. Hence $(\omega, \theta) = (\omega,0) + (0,\theta)\in \iMAforms k$
showing that $((c,b),(h,e),(\omega,\theta))\in \Image (j)$
i.e. $\Ker(p) \subset \Image (r)$. Conversely, if 
$((c,b),(h,e),(\omega,\theta))\in \Image (j)$ then $\omega\in
\rizaMforms {k}$. Thus $\Image(r) \subset \Ker (p)$.

Using the same notation for induced maps, consider the long exact
homology sequence of the above sequence.
\[
\xymatrix{ \dots \ar[r] &  \iMforms {k-1}/ \rizaMforms {k-1}
  \ar[r]^-{\delta} &   \MAintcsg k  \ar[r]^{j} & \MAhsg k \ar[r]^-{p} &
  \iMforms {k}/ \rizaMforms {k} \ar[r] & \dots}
\]

We now claim $\Image(p) = 0$. If $ ((c,b),(h,e),(\omega,\theta))$ is a
cocycle in $\MAccheck k$, then $(\omega, \theta) = (c,b) +
\delta(h,e)$ and so $(\omega,\theta)$ has relative integral
periods. Since $(0,\theta)$ also has relative integral periods we
obtain $\omega\in \rizaMforms {k}$. Hence $p(\cocycle)=0$.

Finally $\Ker(\delta)= \Image(p) = 0$ and we obtain the desired sequence.

The induced map $j$ clearly agrees with $J$ and moreover the boundary
map $\delta\colon\iMforms {k-1}/ \rizaMforms {k-1} \ra \MAintcsg k $
is given by $\delta([\eta]) = ((0,0),(\eta,0),(0,\rho^*\eta))$ and hence
agrees with $\phi$. 
\end{proof}

\section*{Acknowledgements}
The second author is funded by an E.U. Marie-Curie Fellowship and is
grateful to the Institut de Recherche Math\'ematiques Avanc\'ee in
Strasbourg for their hospitality.


\begin{thebibliography}{15}
\bibitem{BottTu}
R. Bott and L. Tu, \emph{Differential Forms in Algebraic Topology}, Graduate Texts in Math. (1986), Springer Verlag.
\bibitem{Brylinski}
J.-L. Brylinski, \emph{Loop Spaces, Characteristic Classes and Geometric Quantization}, Progress in Mathematics, Vol. 107, Birkh\"{a}user, 1993.
\bibitem{CheegerSimons}
J. Cheeger and J. Simons, \emph{Differential characters and geometric
  invariants}, Lecture Notes in Math. 1167 (1985), Springer Verlag,
pp. 50-80.
\bibitem{ChernSimons} 
S.-S. Chern and J. Simons, \emph{Characteristic forms and geometric
  invariants}, Ann. of Math. (2) 99 (1974), 48-69.
\bibitem{HarveyLawson}
R. Harvey and B. Lawson, \emph{Lefschetz-Pontrjagin Duality for
  Differential Characters}, An. Acad. Brasil. Cienc. 73 (2001), no.2, 145-159.
\bibitem{HopkinsSinger}
M. J. Hopkins and I. M. Singer, \emph{Quadratic functions in Geometry,
  Topology, and $M$-theory}, math.AT/0211216. 
\bibitem{Zucchini}
R. Zucchini, \emph{Relative Topological Integrals and Relative
Cheeger-Simons Differential Characters}, J. Geom. Phys. 46 (2003), 355-393. 
\end{thebibliography}
\end{document}